\documentclass[12pt]{article}

\usepackage{amsmath}
\usepackage{amssymb}
\usepackage{eucal}

\usepackage[english]{babel}

\begin{document}
 
\baselineskip 16pt
 
\title{On the Frattini lemma}


\author{V.\,S. Monakhov}


\maketitle

\begin{abstract}

Let $K$ be a subgroup of a finite group $G$, and suppose that $G=KN_G(P)$ 
for every Sylow subgroup $P$ of $K$. Then the subgroup $K$ is normal in 
$G$.

\end{abstract}


MSC2010 20D10, 20D20, 20E34

\medskip

The following assertion is often used in group theory.

\medskip

{\bf Lemma (Frattini).} {\sl
Let $K$ be a normal subgroup of a finite group $G$, and suppose that $P$ is a Sylow
subgroup of $K$. Then $G=KN_G(P)$.}

\medskip

It turns out that the converse of this assertion is also true.

{\bf Lemma.} {\sl
Let $K$ be a subgroup of a finite group $G$, and suppose that $G=KN_G(P)$ 
for every Sylow subgroup $P$ of $K$. Then the subgroup $K$ is normal in 
$G$.}

{\sc Proof.}
Let $P_1,P_2,\ldots P_n$ be the set of all Sylow subgroups of $K$, $x_i\in P_i$, 
$g\in G$. Since by hypothesis, 
$$                                                        
G=KN_G(P_i)=N_G(P_i)K,
$$
for all $i$, it follows that    
$g=a_ib$, $a_i\in N_G(P_i)$, $b\in K$. 
We have $x_i^g=x_i^{a_ib}\in P_i^{a_ib}=P_i^b\subseteq K$.
Since $\langle P_1,P_2,\ldots P_n\rangle =K$,
we see that every element of $K$ can be written as a product of elements 
from $P_1\cup P_2\cup \ldots \cup  P_n$. If $x\in K$ and $g\in G$~are  
arbitrary elements, then 
$$
x=x_1x_2\ldots x_ny_1y_2\ldots y_n\ldots z_1z_2\ldots z_n, \ \ x_i,y_i,z_i\in P_i,
$$
$$
x^g=(x_1x_2\ldots z_n)^g=x_1^gx_2^g\ldots z_n^g\in K,
$$
which means that the subgroup $K$ is normal in $G$. The lemma is proved. 

\bigskip

\noindent V.\,S. MONAKHOV

\noindent Department of mathematics, Gomel F. Scorina State
University, Gomel 246019,  BELARUS 

\noindent E-mail address: Victor.Monakhov@gmail.com

\end{document}